\documentclass[12pt]{article}
\usepackage{mathrsfs}
\usepackage{amsmath}
\usepackage{amsfonts}
\usepackage{amssymb, amsmath, cite}

\newtheorem{theorem}{Theorem}
\newtheorem{lemma}{Lemma}
\newtheorem{proposition}{Proposition}

\newcommand\be{\begin{equation}}
\newcommand\ee{\end{equation}}
\newcommand\ber{\begin{eqnarray}}
\newcommand\eer{\end{eqnarray}}
\newcommand\berr{\begin{eqnarray*}}
\newcommand\eerr{\end{eqnarray*}}

\newcommand{\lm}{\lambda}\newcommand{\bfR}{\mathbb{R}}

\newcommand\re{\mathrm{e}}

\newcommand{\ud}{\mathrm{d}}
\newcommand{\nn}{\nonumber}

\newcommand{\ito}{\int_{\Omega}}

\newcommand{\vep}{\varepsilon}

\title{
Relativistic Chern--Simons--Higgs Vortex Equations}
\author{Xiaosen Han\\Institute of Contemporary Mathematics\\Henan University\\ Kaifeng, Henan 475000, PR China\\ \\Yisong Yang\\Department of Mathematics\\Polytechnic
School of Engineering\\New York University\\ Brooklyn, New York 11201, USA\\ \&\\NYU-ECNU Institute of Mathematical Sciences\\New York University - Shanghai\\3663 North Zhongshan Road, Shanghai 200062, PR China }

\date{}

\begin{document}
\maketitle
\begin{abstract}
An existence theorem is established for the solutions to the non-Abelian relativistic Chern--Simons--Higgs vortex equations over a doubly periodic domain
when the gauge group $G$ assumes the most general and important prototype form, $G=SU(N)$.

\medskip


\end{abstract}

\section{Introduction}
\setcounter{equation}{0}

Let $K=(K_{ij})$ be the Cartan matrix of a semi-simple Lie algebra $L$. Recall that the Toda system is a system of nonlinear elliptic equations over $\bfR^2$, of exponential nonlinearities, of the form
\be \label{I0}
\Delta u_i=-\lm\sum_{j=1}^r K_{ij}\re^{u_j},\quad i=1,\dots,r,
\ee
where $r$ is the rank of $L$. This system is known to be integrable in general \cite{Gan,GGO,Ko,L1,L2,L3,Mi,OT} and arises in the study of non-Abelian monopoles \cite{GGO,OT,WYi}
and nonrelativistic Chern--Simons--Higgs vortices \cite{D0,D,D2,D3,hz}. Interestingly, when $r=1$, that is, when the Cartan subalgebra of $L$ is Abelian such as when
$L$ is the Lie algebra of $SU(2)$, (\ref{I0}) reduces to the classical
Liouville equation \cite{Lio}
\be
\Delta u=-\lm\re^u,
\ee
whose solutions may be constructed by all the integration methods known, such as separation of variables, inverse scattering, the B\"{a}cklund transformation, etc., and is often used as
an illustrative example. It is well known that the Liouville equation and its extensions arise also in differential geometry \cite{aubi,Calabi,KW1,KW2} and have been the focus of various studies
on their analytic aspects \cite{BM,ChangY,CC1,CC2,CC3,CL1,CL2,CN,Li}.
More recently, some extensive work on (\ref{I0}) has been carried out as well aimed at the classification of solutions \cite{Jost,LWY}, understanding its fine analytic structures \cite{ma1,ma2,ma3},
and establishment of bubbling behavior of solutions \cite{JLW,LZhang,OS}.

On the other hand, however, when one considers relativistic Chern--Simons--Higgs vortices \cite{HKP,JW,yang1}, the governing system of equations \cite{D,D2,D3,yang1} is
\be \label{I01}
\Delta u_i=\lm\sum_{j=1}^r\sum_{k=1}^r K_{kj}K_{ji}\re^{u_j}\re^{u_k}-\lm\sum_{j=1}^r K_{ji}\re^{u_j},\quad i=1,\dots,r,
\ee
which deviates from (\ref{I0}) significantly in that, even in the scalar case where $r=1$, the system is nonintegrable \cite{Schiff}. Thus, it is perceivable that the analytic structure of
(\ref{I01}) would be much more complicated than that of (\ref{I0}). Due to its applicability in anyon physics \cite{F,FM1,FM2,Wil}
and challenging mathematical content, it will be desirable to develop an existence theory for the
relativistic vortex equations (\ref{I01}), in which the sources terms resulting from the presence of vortices are temporarily neglected in order to facilitate our discussion.

In \cite{Y0}, a systematic study is conducted to establish an existence theorem for the so-called topological solutions,
realizing the spontaneous symmetry breaking or the celebrated Higgs mechanism of the model, over the full plane $\bfR^2$, to the equations (\ref{I01}), where the Cartan
matrix is of a general form. Another type of solutions of great interest are called the Abrikosov vortices \cite{Ab} or vortex condensates for which the equations are to be solved over a doubly-periodic
lattice domain. In field-theoretic formalism, such a structure is realized by imposing the 't Hooft \cite{tH} boundary condition on gauge and matter fields \cite{SY,WY,yang1}.
Surprisingly, although the system is now considered over a compact domain, double periodicity greatly complicates the problem as evidenced \cite{O} already in the nonrelativistic Liouville equation situation where
one needs to use the Weierstrass elliptic functions as building blocks for solutions. In the relativistic situation, although there is a variational principle, the action functional is not bounded
from below and one has to consider a constrained minimization problem. However, the presence of the constraints leads to a Lagrange multiplier issue so that it prevents one from recovering
the original equations of motion. In order to overcome this difficulty, one has to consider an inequality-constrained problem instead and bypass the Lagrange multiplier problem with achieving
an interior minimum. As a consequence, the progress in developing an existence theory for the doubly periodic solutions of the relativistic Chern--Simons--Higgs vortex equations has been
slow and sporadic. Specifically, in \cite{caya1}, an inequality-constrained minimization method was first used to establish the existence of solutions for the scalar case of (\ref{I01}), namely
$r=1$ or $G=SU(2)$, and, in \cite{nota}, the method in \cite{caya1} was remarkably extended and refined to tackle the first non-scalar case of (\ref{I01}), namely $r=2$ or $G=SU(3)$, in which
there are two inequality constraints characterizing the solvability of two quadratic constraints. The study on the next important non-scalar case of (\ref{I01}), that is, $r=3$ or $G=SU(4)$,
had been unsuccessful due to the difficulty in resolving more than two quadratic constraints simultaneously until
very recently a new idea based on an implicit-function theorem argument was implemented to resolve three coupled quadratic constraints, which enables the establishment of an existence theorem for
$G=SU(4)$ in \cite{HLY}.

The purpose of the present paper is to establish an existence theorem for the doubly periodic solutions of (\ref{I01}) for the most general situation, $G=SU(N)$ ($N\geq2$) or $r=N-1\geq1$,
adapting the implicit function method initiated in \cite{HLY} in handling multiple quadratic constraints. Instead of a `squeeze-to-the-middle' approach in \cite{HLY} for
resolving the constraints, however, we use here an ordered iterative scheme which is effective and easier to implement in the general situation.
In the next section, we state our main existence theorem. In the section after, we prove the theorem.
In the last section we end the paper with some remarks.

{\em Note added upon acceptance.} This work grew out of an earlier version of \cite{HLY} in which only the $SU(4)$ problem was resolved. Since the submission of the present
work, further development has been achieved to resolve the general situation when $G$ is a simple Lie group by using a degree-theory method formalism, which 
greatly expands
\cite{HLY} into its current
 updated version seen.

\section{Existence theorem  }
\setcounter{equation}{0}
\setcounter{lemma}{0}

With the source terms in the presence of multiply distributed vortices, the non-Abelian relativistic Chern--Simons--Higgs equations are \cite{D,D2,D3,yang1}
\be \label{I1}
\Delta u_i=\lm\left(\sum_{j=1}^r\sum_{k=1}^r K_{kj}K_{ji}\re^{u_j}\re^{u_k}-\sum_{j=1}^r K_{ji}\re^{u_j}\right)+4\pi\sum_{j=1}^{N_i}\delta_{p_{ij}}(x),\quad i=1,\dots,r,
\ee
where $\delta_p$ denotes the Dirac measure concentrated at the point $p$, $\lm>0$ is a coupling constant, and
the equations are considered over a doubly periodic domain $\Omega$ resembling a lattice cell housing a distribution of dually charged vortices located at $p_{ij},j=1,\dots,N_i, i=1,\dots,r$.

  We focus on the  system \eqref{I1}  with  $G=SU(n+1),  n\ge2$.   The Cartan matrix $K$ for $SU(n+1)$ is an $n\times n$ matrix  given by
  \ber
 K\equiv\begin{pmatrix}
 2&-1&0 &\cdots&0\\
 -1&2&-1 &\cdots&0\\
  \vdots&\ddots&\ddots&\ddots&\vdots\\
  0&\dots&-1&2&-1\\
 0& \dots& 0&-1&2
\end{pmatrix}.\label{t2}
\eer

It is easy to check that  $K^{-1}$ is symmetric with entries given by
 \be
 \quad (K^{-1})_{ij}=\frac{i(n+1-j)}{n+1}, \quad i\le j,\,  i=1, \dots, n. \label{s10}
 \ee

Then after the translation
\be
u_i\to u_i+\ln R_i \quad \text{with} \quad R_i\equiv \sum\limits_{j=1}^n(K^{-1})_{ij}=\frac{i(n+1-i)}2,\quad i=1, \dots, n, \label{t4}
\ee
 the system  \eqref{I1} can be rewritten as
 \be
\Delta u_i=\lm\left(\sum_{j=1}^n\sum_{k=1}^n \tilde{K}_{jk}\tilde{K}_{ij}\re^{u_j}\re^{u_k}-\sum_{j=1}^n \tilde{K}_{ij}\re^{u_j}\right)
+4\pi\sum_{j=1}^{N_i}\delta_{p_{ij}}(x),\quad i=1,\dots,n,\label{t4a}
\ee
or in a vector form,
\ber
 \Delta \mathbf{u}=\lambda\tilde{K}\mathrm{U}\tilde{K}(\mathbf{U}-\mathbf{1}) +4\pi\mathbf{s}, \label{t5}
\eer
 where

 \ber
&&\tilde{K}\equiv KR= \begin{pmatrix}
 n&-(n-1)&0 &\cdots&\cdots&0\\
 -\frac{n}{2}&2(n-1)& -\frac{3(n-2)}{2} &\cdots&\cdots&0\\
  \vdots&\ddots&\ddots&\ddots&\vdots&\vdots\\
   \dots&\dots&-\frac{(i-1)(n+2-i)}{2} &i(n+1-i)&-\frac{(i+1)(n-i)}{2}&\cdots\\
  \vdots&\vdots&\ddots&\ddots&\ddots&\vdots\\
  0&\dots&\cdots&-\frac{3(n-2)}{2} &2(n-1)&-\frac{n}{2}\\
 0& \dots& \cdots&0&-(n-1)&n
\end{pmatrix}, \label{t7}\eer
\ber
&&R\equiv{\rm diag}\left\{R_1, \dots,  R_n\right\},  \quad  \mathbf{1}=(1, \dots, 1)^\tau,\label{t8}\\
  &&\mathbf{u}=(u_1, \dots, u_n)^\tau, \quad  \mathrm{U}={\rm diag}\left\{\re^{u_1},  \dots, \re^{u_n}\right\},\quad \mathbf{U}=(\re^{u_1},  \dots, \re^{u_n})^\tau,  \label{t6}\\
&&\mathbf{s}=\left(\sum\limits_{s=1}^{N_1}\delta_{p_{1s}},  \dots,   \sum\limits_{s=1}^{N_n}\delta_{p_{ns}}\right)^\tau.\label{t9}
\eer

By the definition of $\tilde{K}$ given in \eqref{t7},  we  obtain the following simple facts
 \ber
  \quad \tilde{K}^{-1}=R^{-1}K^{-1}, \quad  \tilde{K}^{-1}\mathbf{1}=\mathbf{1},\quad K^{-1}\mathbf{1}=R\mathbf{1},\label{t6a}
 \eer
which will be repeatedly  used  later in this paper.

We are interested in the existence of solutions of  \eqref{t4a} or  \eqref{t5}   over a doubly periodic domain $\Omega$.
Our main result  reads as follows.

\begin{theorem}\label{th2}
Consider the  nonlinear  elliptic system \eqref{t4a} or  \eqref{t5} over   a doubly periodic domain $\Omega$ in $\mathbb{R}^2$. For any  given   points $p_{i1}, \dots,  p_{iN_i} \in \Omega\, (i=1, \dots n)$,
 which need not to  be distinct, the following conclusions hold.

 {\rm (i)} (Necessary condition for existence) If
    \be
     \lambda\le\lambda_0\equiv\frac{16\pi}{|\Omega|}\frac{\sum\limits_{i=1}^n\sum\limits_{j=1}^n(K^{-1})_{ij}N_j}{\sum\limits_{i=1}^n\sum\limits_{j=1}^n(K^{-1})_{ij}},  \label{t10}
    \ee
there is no solution   to the system, where the entries of  $K^{-1}$ are given by \eqref{s10}. In other words, a solution can exist only when $\lm$ is larger than the right-hand side of
 \eqref{t10}.

{\rm (ii)} (Sufficient condition for existence)   There exists some    $\lambda_1>\lambda_0$  such that when $\lambda>\lambda_1$ the  system  admits a  solution  over $\Omega$.

 {\rm (iii)} (Asymptotic behavior)  The solution $(u_1, \dots, u_n)$ obtained above satisfies
   \be
    \lim\limits_{\lambda\to\infty}\ito \left(\re^{u_i}-1\right)^2\ud x=0, \quad i=1, \dots, n.
    \ee

 {\rm (iv)} (Quantized integrals)    If $(u_1, \dots, u_n)$ is a solution then there hold the quantized integrals
  \ber
   \ito\left(\sum_{j=1}^n\sum_{k=1}^n \tilde{K}_{jk}\tilde{K}_{ij}\re^{u_j}\re^{u_k}-\sum_{j=1}^n \tilde{K}_{ij}\re^{u_j}\right)\ud x=-\frac{4\pi N_i}{\lambda}, \quad i=1, \dots, n.\label{t11}
  \eer
\end{theorem}

 Note that when the vortex numbers  $N_i\,(i=1, \dots, n)$ are the same, say $m$, our necessary condition \eqref{t10} reduces to that for the $U(1)$ case as in \cite{caya1} as follows
  \be
   \lambda\le\frac{16\pi m}{|\Omega|}.
  \ee
which is a comfort and also surprising since we are now considering a non-Abelian and non-scalar situation.

\section{Proof of theorem}
 \setcounter{equation}{0}
  \setcounter{lemma}{0}

 In this section  we   apply  a constrained minimization  procedure  developed in \cite{caya1}, which was later  modified in  \cite{nota},  to  establish the existence of doubly periodic solutions to
(\ref{t4a})  \eqref{t5} when $n=2$. We carry out our proof in several steps. First we show that the condition \eqref{t10} implies nonexistence of solutions as stated and we explore a variational structure of our
equations. Next we show how to resolve multiple constraints using an iterative scheme and an implicit function argument.
We then conduct a constrained minimization procedure and show that there is a solution to the minimization problem. In the subsequent subsection, we establish some suitable estimates
which ensures that the minimum point obtained must be an interior minimum, thus ruling out the Lagrange multiplier issue.
In the last two subsections, we show that the interior minimum point obtained is a classical solution of the original vortex equations and we then establish the stated asymptotic behavior
of solutions as $\lm\to\infty$ and quantized integrals.

\subsection{Necessary condition and variational structure}
Let $u_i^0$ be the solution of the following problem (see \cite{aubi})
\ber
 \Delta u_i^0=4\pi\sum\limits_{s=1}^{N_i}\delta_{p_{is}}-\frac{4\pi N_i}{|\Omega|},\quad \ito u_i^0\ud x=0, \quad i=1, \dots,  n,
\eer
and $u_i=u_i^0+v_i, \, i=1,  \dots,  n$.  Let us introduce the notation ($n$-vectors)
\be
 \mathbf{v}=(v_1, \dots, v_n)^\tau, \quad \mathbf{N}= (N_1, \dots, N_n)^\tau, \quad \mathbf{0}=(0, \dots, 0)^\tau. \label{s2}
 \ee
Then   we   reformulate  the  system  \eqref{t4a} or \eqref{t5} as
 \be
\Delta v_i=\lm\left(\sum_{j=1}^n\sum_{k=1}^n \tilde{K}_{jk}\tilde{K}_{ij}\re^{u_j^0+v_j}\re^{u_k^0+v_k}-\sum_{j=1}^n \tilde{K}_{ij}\re^{u_j^0+v_j}\right)+\frac{4\pi N_i}{|\Omega|},
\quad i=1,\dots,n,\label{s2a}
\ee
or
\ber
 \Delta \mathbf{v}=\lambda\tilde{K}\mathrm{U}\tilde{K}(\mathbf{U}-\mathbf{1})+ \frac{4\pi}{|\Omega|} \mathbf{N},  \label{s1}
\eer
with the understanding  that
\ber
\mathrm{U}={\rm diag}\left\{\re^{u_1^0+v_1},  \dots, \re^{u_n^0+v_n}\right\}, \quad \mathbf{U}=(\re^{u_1^0+v_1}, \dots, \re^{u_n^0+v_n})^\tau. \label{s1'}
\eer

 We first present   a  necessary  condition for the existence of solution to \eqref{s2a} or  \eqref{s1}.
   For any solution $ \mathbf{v} $ of  \eqref{s1},   taking integration over $\Omega$, we obtain the natural constraint
 \ber
 \ito\tilde{K}\mathrm{U}\tilde{K}(\mathbf{U}-\mathbf{1})\ud x+\frac{4\pi}{\lambda}\mathbf{N}=\mathbf{0}. \label{s11a}
 \eer
 Multiplying both sides of \eqref{s11a} by $K^{-1}$, we have
 \ber
  \ito R\mathrm{U}KR(\mathbf{U}-\mathbf{1})\ud x+\frac{4\pi}{\lambda}K^{-1}\mathbf{N}=\mathbf{0}.\label{s12}
\eer
 Then  noting \eqref{s12}, the positive definiteness of $K$, and the fact that $K^{-1}\mathbf{1}=R\mathbf{1}$, we have
  \ber
 0&=&\ito (R\mathbf{U})^\tau KR(\mathbf{U}-\mathbf{1})\ud x+\frac{4\pi}{\lambda}\mathbf{1}^\tau K^{-1}\mathbf{N}\nn\\
  &=& \ito \left(R\left[\mathbf{U}-\frac12{\mathbf{1}}\right]\right)^\tau K\left(R\left[\mathbf{U}-\frac12{\mathbf{1}}\right]\right)\ud x-\frac{|\Omega|}{4}(R\mathbf{1})^\tau K(R\mathbf{1})+\frac{4\pi}{\lambda}\mathbf{1}^\tau K^{-1}\mathbf{N}\nn\\
   &>&-\frac{|\Omega|}{4}\mathbf{1}^\tau K^{-1}\mathbf{1} +\frac{4\pi}{\lambda}\mathbf{1}^\tau K^{-1}\mathbf{N}, \label{s13}
\eer
which implies
 \ber
  \lambda>  \frac{16\pi \mathbf{1}^\tau K^{-1}\mathbf{N}}{|\Omega|\mathbf{1}^\tau K^{-1}\mathbf{1} }. \label{s18}
 \eer
That is,  \eqref{s18} spells out a necessary condition for the existence of solutions to \eqref{s1}. Hence the first conclusion  of Theorem \ref{th2} follows.

Now we  aim to  find a variational principle for the equations  \eqref{s1}.  After a simple computation, we see  that the matrix $\tilde{K}$   admits  a decomposition
  \ber
   \tilde{K}=\tilde{P}\tilde{S},\label{s3}
  \eer
 with
  \ber
  \tilde{P}\equiv{\rm diag}\left\{\tilde{P}_1,\dots, \tilde{P}_n\right\}, \quad \tilde{P}_i\equiv \frac{n}{i(n+1-i)}, \quad i=1, \dots, n, \label{s4}
  \eer
  \ber
  \tilde{S} \equiv \begin{pmatrix}
 n&-(n-1)&0 &\cdots&0\\
  -(n-1)&\frac{4(n-1)^2}{n}&-\frac{3(n-1)(n-2)}{n}&\dots&\dots\\
 \quad \ddots&\quad\ddots&\quad\ddots&\vdots&\vdots\\
   \dots&-\tilde{S}_{ii-1} &\tilde{S}_{ii}&-\tilde{S}_{ii+1}&\cdots\\
\vdots&\ddots&\ddots&\ddots&\vdots\\
\dots&\dots&-\frac{3(n-1)(n-2)}{n}&\frac{4(n-1)^2}{n}& -(n-1)\\
 0& \cdots&0& -(n-1)&n
\end{pmatrix},   \label{s5}
  \eer
which is a tridiagonal matrix with
 \be\label{s6}
\left\{\begin{array}{rcl}
  \tilde{S}_{ii-1}&\equiv&\frac{(i-1)i(n+2-i)(n+1-i)}{2n},\\
 \tilde{S}_{ii}&\equiv&\frac{i^2(n+1-i)^2}{n}, \\
 \tilde{S}_{ii+1}&\equiv&\frac{i(i+1)(n+1-i)(n-i)}{2n},\end{array}\right. i=3, \quad\dots, n-2.
 \ee

  Then we may rewrite \eqref{s1}  equivalently  as
 \ber
 \Delta M\mathbf{v}&=&\lambda\mathrm{U}\tilde{S}(\mathbf{U}-\mathbf{1})+\frac{\mathbf{b}}{|\Omega|}, \label{s8}
\eer
where
 \ber
M&\equiv& \tilde{P}^{-1}\tilde{K}^{-1}=\frac{2}{n}K^{-1}, \quad \tilde{P}^{-1}= {\rm diag}\left\{\tilde{P}_1^{-1}, \dots, \tilde{P}_n^{-1}\right\}, \label{s9}\\
   \tilde{P}_i^{-1}&=&\frac{i(n+1-i)}{n},\quad i=1, \dots, n, \label{s9a}
 \eer
 and
 \be
 \mathbf{b}=(b_1, \dots, b_n)^\tau\equiv4\pi M\mathbf{N}.\label{s10a}
  \ee
By the definition \eqref{s10a} for $ \mathbf{b}$,  we easily   find  that    $b_i>0,\, i=1,\dots, n$.

 We use $W^{1, 2}(\Omega)$ to denote the  Sobolev space of scalar-valued or vector-valued $\Omega$-periodic $L^2$ functions
  with their derivatives also belonging to   $L^2(\Omega)$.

It may be examined that the equations   \eqref{s8} are the Euler--Larange equations of the action functional
 \ber
 I(\mathbf{v})&=& \frac12\sum\limits_{i=1}^2\ito \partial_i \mathbf{v}^\tau M\partial_i\mathbf{v}\ud x
 +\frac\lambda2\ito (\mathbf{U}-\mathbf{1})^\tau \tilde{S}(\mathbf{U}-\mathbf{1})\ud x+ \frac{1}{|\Omega|}\ito \mathbf{b}^\tau\mathbf{v}\ud x,\label{s11}
\eer
where  we use the notation  \eqref{s5}, \eqref{s9}  and  \eqref{s10a} throughout this paper.

In the following subsections  we will   use  a constrained minimization approach to find  the   critical points of the functional $I$.

\subsection{Multiple constraints}
To start our constrained minimization process, we need to find some suitable constraints subject to which the functional $I$ will be minimized.

  Note   the space   $W^{1,2}(\Omega)$   can be decomposed as
  \be
  W^{1,2}(\Omega)=\mathbb{R}\oplus  \dot{W}^{1,2}(\Omega),
  \ee
  where
\be \dot{W}^{1,2}(\Omega)\equiv\left\{w\in W^{1,2}(\Omega)\Bigg| \ito w\ud x=0\right\},
\ee
  is a closed subspace of $W^{1,2}(\Omega)$.

Then, for  $v_i\in W^{1,2}(\Omega)$,  we  have the decomposition
 \be \label{**}
v_i=c_i+w_i,\quad   c_i \in \mathbb{R}, \quad w_i\in  \dot{W}^{1,2}(\Omega),\quad i=1, \dots, n.
\ee

If $\mathbf{v}\in W^{1,2}(\Omega)$ satisfies the constraint  \eqref{s11a}, which is equivalent to
 \be
  \ito \mathrm{U}\tilde{S}(\mathbf{U}-\mathbf{1})\ud x+ \frac{\mathbf{b}}{\lambda}=\mathbf{0}, \label{as1}
 \ee
then by the  decomposition (\ref{**}) with ${\bf w}=(w_1,\dots,w_n)^\tau$,  we obtain
 \ber
  &&\re^{2c_1}a_{11}-\re^{c_1}P_1(\mathbf{w}; \re^{c_2})+\frac{b_1}{n\lambda}=0,\label{s19}\\
  &&\re^{2c_i}a_{ii}-\re^{c_i}P_i(\mathbf{w}; \re^{c_{i-1}},  \re^{c_{i+1}})+\frac{nb_i}{i^2(n+1-i)^2\lambda}=0,\quad i=2, \dots, n-1,\label{s20}\\
  &&\re^{2c_n}a_{nn}-\re^{c_n}P_n(\mathbf{w}; \re^{c_{n-1}})+\frac{b_n}{n\lambda}=0,\label{s21}
 \eer
where  and in the sequel the notation
\ber
P_1(\mathbf{w}; \re^{c_2})&\equiv& \frac{a_1}{n}+\frac{(n-1)a_{12}}{n}\re^{c_2},\label{s22}\\
P_i(\mathbf{w}; \re^{c_{i-1}}, \re^{c_{i+1}})&\equiv& \frac{a_i}{i(n+1-i)}+ \frac{(i-1)(n+2-i)a_{ii-1}}{2i(n+1-i)}\re^{c_{i-1}} \nn\\
&&+ \frac{(i+1)(n-i)a_{ii+1}}{2i(n+1-i)}\re^{c_{i+1}},\quad i=2, \dots, n-1,\label{s23}\\
P_n(\mathbf{w}; \re^{c_{n-1}})&\equiv& \frac{a_n}{n}+\frac{(n-1)a_{nn-1}}{n}\re^{c_{n-1}},\label{s24}
\eer
and the definition
\ber
a_i&\equiv& a_i(w_i)\equiv\ito\re^{u_i^0+w_i}\ud x,\label{s24a}\\ a_{ij}&\equiv& a_{ij}(w_i, w_j)\equiv\ito\re^{u_i^0+u_j^0+w_i+w_j}\ud x,
  \quad i, j=1, \dots, n, \label{s24b}
\eer
are used.

Then, we see that, for any $\mathbf{w}\in \dot{W}^{1,2}(\Omega)$, the equations \eqref{s19}--\eqref{s21} are solvable in
\be
\mathbf{c}\equiv(c_1, \dots, c_n)^\tau,
\ee   only if
 \ber
   P_1^2(\mathbf{w}; \re^{c_2})&\ge& \frac{4b_1a_{11}}{n\lambda}, \label{s25}\\
   P_i^2(\mathbf{w};  \re^{c_{i-1}}, \re^{c_{i+1}})&\ge&\frac{4nb_ia_{ii}}{i^2(n+1-i)^2\lambda}, \quad  i=2, \dots, n-1, \label{s26}\\
    P_n^2(\mathbf{w}; \re^{c_{n-1}})&\ge& \frac{4b_na_{nn}}{n\lambda},  \label{s27}
 \eer
are satisfied which can be ensured by the following  simpler inequality-type constraints
 \ber
  a_i^2&\ge& \frac{4nb_ia_{ii}}{\lambda}, \quad i=1,  \dots, n,\label{s28}
 \eer
or expressed explicitly as
 \ber
   \left(\ito\re^{u_i^0+w_i}\ud x\right)^2&\ge& \frac{4nb_i}{\lambda}\ito\re^{2u_i^0+2w_i}\ud x, \quad i=1,  \dots, n.\label{s28'}
 \eer

Now we define   the  admissible  set
 \be
  \mathcal{A}\equiv\Big\{\mathbf{w}\Big|\, \mathbf{w}\in \dot{W}^{1, 2}(\Omega)\, \,\text{satisfies}\,\, \eqref{s28'}\Big\}. \label{s29}
 \ee

Thus,    for any  $\mathbf{w}\in\mathcal{A}$,  we can find a solution of   the equations \eqref{s19}--\eqref{s21} with respect to $(c_1,  \dots,  c_n)$   by solving    the following  coupled equations
 \ber
 \re^{c_1}&=&\frac{ P_1(\mathbf{w};  \re^{c_2})+\sqrt{P_1^2(\mathbf{w};  \re^{c_2})-\frac{4b_1a_{11}}{n\lambda}}}{2a_{11}}
 \nn\\&\equiv& f_1(\re^{c_2}), \label{s30}\\
  \re^{c_i}&=&\frac{ P_i(\mathbf{w};  \re^{c_{i-1}}, \re^{c_{i+1}})+\sqrt{P^2_i(\mathbf{w};  \re^{c_{i-1}}, \re^{c_{i+1}})-\frac{4nb_ia_{ii}}{i^2(n+1-i)^2\lambda}}}{2a_{ii}}
 \nn\\&\equiv& f_i(\re^{c_{i-1}}, \, \re^{c_{i+1}}),\quad i=2, \dots, n-1,\label{s31}\\
  \re^{c_n}&=&\frac{ P_n(\mathbf{w};  \re^{c_{n-1}})+\sqrt{P_n^2(\mathbf{w};  \re^{c_{n-1}})-\frac{4b_na_{nn}}{n\lambda}}}{2a_{nn}}
 \nn\\&\equiv& f_n(\re^{c_{n-1}}).\label{s32}
 \eer

For any $\mathbf{w}\in  \mathcal{A}$, to  solve the system of the equations  \eqref{s30}--\eqref{s32} with respect to $(c_1, \dots, c_n)$,  we will convert the system into a single equation.

 For $n\ge 2$,   we need to find a solution of the system
 \ber
 &&t_1-f_1(t_2)=0,\label{r1}\\
&& t_i-f_i(t_{i-1}, t_{i+1})=0,\quad  i=2, \dots, n-1,\label{r2}\\
 &&t_n-f_n(t_{n-1})=0.  \label{r3}
 \eer

  \begin{proposition}\label{prop1}
   For any $ \mathbf{w} \in   \mathcal{A}$ and $n\ge2$, the system \eqref{r1}--\eqref{r3} admits a unique solution in $(0, \infty)^n$.
 \end{proposition}

  {\bf Proof. }  When $n=2, 3$,  the system can be transformed into a single equation directly, which was solved in  \cite{nota} and \cite{HLY}, respectively.
However, for $n\ge4$,  it difficult to reduce the system \eqref{r1}--\eqref{r3} into a single equation directly.  Here we  use the implicit function theorem to  overcome this difficulty.

  For $n\ge 4$, using the first two equations of the system  \eqref{r1}--\eqref{r3}, we have  the   relation
 \ber
  F_2(t_2, t_3)\equiv t_2-f_2(f_1(t_2), t_3)=0, \quad t_2, t_3>0.  \label{ar7}
 \eer

   We first show that the relation \eqref{ar7} may uniquely determine an implicit function
  \be
   t_2=g_2(t_3)>0,\quad t_3>0. \label{ar6}
  \ee

By the expression \eqref{s30}--\eqref{s32}, we have
 \ber
  \frac{\ud f_1(t_2)}{\ud t_2}&=&\frac{1}{2a_{11}}\left\{ \frac{(n-1)a_{12}}{n}+\frac{\frac{(n-1)a_{12}}{n}P_1(\mathbf{w}; t_2)}{\sqrt{P_1^2(\mathbf{w};  t_2)-\frac{4b_1a_{11}}{n\lambda}}} \right\}\nn\\
  &=&\frac{\frac{(n-1)a_{12}}{n}f_1(t_2)}{\sqrt{P_1^2(\mathbf{w};  t_2)-\frac{4b_1a_{11}}{n\lambda}}}. \label{r7}
 \eer
 Similarly, we obtain
  \ber
   \frac{\partial f_2(t_1, t_3)}{\partial t_1}=\frac{\frac{na_{21}}{4(n-1)}f_2(t_1, t_3)}{\sqrt{P_2^2(\mathbf{w};  t_1, t_3)-\frac{nb_2a_{22}}{(n-1)^2\lambda}}}. \label{r8}
  \eer

Then, using \eqref{r7}--\eqref{r8} and the constraints \eqref{s28}, we  have
 \ber
  \frac{\partial F_2(t_2, t_3)}{\partial t_2}&=&1-\frac{\partial f_2(t_1, t_3)}{\partial t_1}\Big|_{t_1=f_1(t_2)} \frac{\ud f_1(t_2)}{\ud t_2}\nn\\
   &=&1-\frac{\frac14a_{12}^2f_1(t_2)f_2(f_1(t_2), t_3)}{\sqrt{P_1^2(\mathbf{w};  t_2)-\frac{4b_1a_{11}}{n\lambda}}\sqrt{P_2^2(\mathbf{w};  f_1(t_2), t_3)-\frac{nb_2a_{22}}{(n-1)^2\lambda}} }\nn\\
   &>& 1-\frac{\frac14a_{12}^2f_1(t_2)f_2(f_1(t_2), t_3)}{\frac{{(n-1)}a_{12}}{n}t_2\frac{na_{12}}{4(n-1)}f_1(t_2)}\nn\\
   &=&1-\frac{f_2(f_1(t_2), t_3)}{t_2}=\frac{F_2(t_2, t_3)}{t_2}=0, \quad t_2, t_3>0.\label{r9}
 \eer

Therefore,  from \eqref{r9}  and the implicit function theorem we see that there exists a unique implicit function
\ber
 t_2=g_2(t_3)>0,\quad t_3>0,\quad  \text{such that }\, F_2(g_2(t_3), t_3)\equiv0,  \quad t_3>0. \label{r10}
  \eer

Similarly, combining  \eqref{r10} and the third $(i=3)$ equation of the system \eqref{r1}--\eqref{r3} gives the relation
 \ber
  F_3(t_3, t_4)\equiv t_3-f_3(g_2(t_3), t_4)=0, \quad t_3, t_4>0,
 \eer
which  also determines a unique implicit function
 \be
  t_3=g_3(t_4)>0,\quad t_4>0,  \label{ar9}
 \ee
  such that
  \be
   F_3(g_3(t_4), t_4)\equiv0,  \quad t_4>0. \label{ar9'}
  \ee

Repeating the above   procedure,  we obtain that there exists a family of uniquely determined implicit functions
 \ber
   t_i=g_i(t_{i+1})>0,\quad t_{i+1}>0, \quad i=4, \dots, n-1, \label{ar10}
 \eer
satisfying
 \ber
   F_i(t_i, t_{i+1})\equiv t_i-f_i(g_{i-1}(t_i), t_{i+1})=0, \quad t_i,\,  t_{i+1}>0,  \quad i=4, \dots, n-1. \label{ar11}
 \eer

Therefore,  in view of \eqref{r10}, \eqref{ar9}--\eqref{ar11},  to solve the system \eqref{r1}--\eqref{r3}, it is equivalent to solve the following single equation
 \ber
 F(t_n) \equiv t_n-f_n\left(g_{n-1}(t_n)\right)=0, \quad  t_n\ge0. \label{ar14}
 \eer

In other words, to prove our proposition,  we just need to  show that $F(\cdot)$ admits a unique positive zero.

We easily see that
    \be
     F(0)<0. \label{ar15}
    \ee

Next we prove that
  \ber
  \lim\limits_{t_n\to \infty}F(t_n)=\infty. \label{ar16}
 \eer

 Noting \eqref{r1}--\eqref{r3}, after a direct computation,   we obtain
  \ber
   \lim\limits_{t_n\to\infty}\frac{t_1}{t_n}&=&\frac{(n-1)a_{12}}{na_{11}} \lim\limits_{t_n\to\infty}\frac{t_2}{t_n},\label{rr1}\\
    \lim\limits_{t_n\to\infty}\frac{t_i}{t_n}&=&\frac{1}{a_{ii}} \left\{\frac{(i-1)(n+2-i)a_{ii-1}}{2i(n+1-i)}\lim\limits_{t_n\to\infty}\frac{t_{i-1}}{t_n}\right.\nn\\
    &&\left.+ \frac{(i+1)(n-i)a_{ii+1}}{2i(n+1-i)}\lim\limits_{t_n\to\infty}\frac{t_{i+1}}{t_n} \right\}, \quad i=2, \dots, n-2,\label{rr2}\\
     \lim\limits_{t_n\to\infty}\frac{t_{n-1}}{t_n}&=&\frac{1}{a_{n-1n-1}} \left\{\frac{3(n-2)a_{n-1n-2}}{4(n-1)}\lim\limits_{t_n\to\infty}\frac{t_{n-2}}{t_n}+\frac{na_{n-1n}}{4(n-1)}\right\}.\label{rr3}
  \eer

 Then  we infer from  \eqref{rr1}--\eqref{rr3} and H\"{o}lder's inequality that
  \ber
   \lim\limits_{t_n\to\infty}\frac{t_i}{t_n}&\le& \frac{(n-i)a_{ii+1}}{(n-i+1)a_{ii}}\lim\limits_{t_n\to\infty}\frac{t_{i+1}}{t_n},\quad i=2,\dots, n-1.\label{rr4}
  \eer

 In particular,  we have
  \ber
  \lim\limits_{t_n\to\infty}\frac{g_{n-1}(t_n)}{t_n}=\lim\limits_{t_n\to\infty}\frac{t_{n-1}}{t_n}\le \frac{a_{n-1n}}{2a_{n-1n-1}}.\label{rr5}
  \eer

Then  using  \eqref{rr5}  and H\"{o}lder's inequality again, we have
 \ber
  \lim\limits_{t_n\to\infty}\frac{F(t_n)}{t_n}&=&1-\lim\limits_{t_n\to\infty}\frac{f_n(g_{n-1}(t_n))}{t_n}\nn\\
    &=&1-\frac{(n-1)a_{nn-1}}{na_{nn}}\lim\limits_{t_n\to\infty}\frac{g_{n-1}(t_n)}{t_n}\nn\\
     &\ge& 1- \frac{(n-1)a^2_{nn-1}}{2na_{n-1n-1}a_{nn}}\nn\\
     &\ge&1-\frac{n-1}{2n}\nn\\
      &=&\frac{n+1}{2n}. \label{ar19}
 \eer

Hence \eqref{ar19} implies the desired limit  \eqref{ar16}.

 Consequently,  from  \eqref{ar15}--\eqref{ar16}, we conclude that $F(\cdot)$ admits at least one positive zero.  In the following, we  prove  the uniqueness of the
 zero.

 Using \eqref{s30}--\eqref{s32} and the constraints \eqref{s28} and after a direct computation,  we have
 \ber
\frac{\partial f_{n-1}(t_{n-2}, t_n)}{\partial t_n}&=&\frac{\frac{na_{n-1n}}{4(n-1)}f_{n-1}(t_{n-2}, t_n)}{\sqrt{P_{n-1}^2(\mathbf{w};  t_{n-2}, t_n)-\frac{nb_{n-1}a_{n-1n-1}}{(n-1)^2\lambda}}}\nn\\
 &\le& \frac{\frac{na_{n-1n}}{4(n-1)}f_{n-1}(t_{n-2}, t_n)}{\frac{3(n-2)a_{n-1n-2}}{4(n-1)}t_{n-2}+\frac{na_{n-1n}}{4(n-1)}t_n},\label{rr8}
\eer
\ber
\frac{\partial f_{n-1}(t_{n-2}, t_n)}{\partial t_{n-2}}&=&\frac{\frac{3(n-2)a_{n-1n-2}}{4(n-1)}f_{n-1}(t_{n-2}, t_n)}{\sqrt{P_{n-1}^2(\mathbf{w};  t_{n-2}, t_n)-\frac{nb_{n-1}a_{n-1n-1}}{(n-1)^2\lambda}}}\nn\\
 &\le& \frac{\frac{3(n-2)a_{n-1n-2}}{4(n-1)}f_{n-1}(t_{n-2}, t_n)}{\frac{3(n-2)a_{n-1n-2}}{4(n-1)}t_{n-2}+\frac{na_{n-1n}}{4(n-1)}t_n},\label{rr9}
\eer
\ber
\frac{\partial f_{n-2}(t_{n-3}, t_{n-1})}{\partial t_{n-1}}&=&\frac{\frac{(n-1)a_{n-2n-1}}{3(n-2)}f_{n-2}(t_{n-3}, t_{n-1})}{\sqrt{P_{n-2}^2(\mathbf{w};  t_{n-3}, t_{n-1})-\frac{4nb_{n-2}a_{n-2n-2}}{9(n-2)^2\lambda}}}\nn\\
  &<& \frac{\frac{(n-1)a_{n-2n-1}}{3(n-2)}f_{n-2}(t_{n-3}, t_{n-1})}{\frac{(n-1)a_{n-2n-1}}{3(n-2)}t_{n-1}}\nn\\
  &=&\frac{f_{n-2}(t_{n-3}, t_{n-1})}{t_{n-1}}.\label{rr10}
\eer

Then, from \eqref{rr8}--\eqref{rr10}, we infer that
\ber
 \frac{\ud g_{n-1}(t_n)}{\ud t_n}&=&\frac{\frac{\partial f_{n-1}(t_{n-2}, t_n)}{\partial t_n}}{1-\frac{\partial f_{n-1}(t_{n-2}, t_n)}{\partial t_{n-2}}\frac{\partial t_{n-2}}{\partial t_{n-1}}}\nn\\
  &<& \frac{\frac{\frac{na_{n-1n}}{4(n-1)}f_{n-1}(t_{n-2}, t_n)}{\frac{3(n-2)a_{n-1n-2}}{4(n-1)}t_{n-2}+\frac{na_{n-1n}}{4(n-1)}t_n}}
  {1-\frac{\frac{3(n-2)a_{n-1n-2}}{4(n-1)}f_{n-1}(t_{n-2}, t_n)}{\frac{3(n-2)a_{n-1n-2}}{4(n-1)}t_{n-2}+\frac{na_{n-1n}}{4(n-1)}t_n}
  \frac{f_{n-2}(t_{n-3}, t_{n-1})}{t_{n-1}}}\nn\\
  &=&\frac{f_{n-1}(t_{n-2}, t_n)}{t_n}=\frac{t_{n-1}}{t_n}=\frac{g_{n-1}(t_n)}{t_n}.\label{rr11}
\eer

Hence, using \eqref{rr11}  and the  constraints \eqref{s28}, we have
 \ber
  \frac{\ud F(t_n)}{\ud t_n}&=&1-\frac{\partial f_n(t_{n-1})}{\partial t_{n-1}}\frac{\ud  t_{n-1}}{\ud  t_n}\Big|_{t_{n-1}=g_{n-1}(t_n)}\nn\\
  &=&1-\frac{\frac{(n-1)a_{nn-1}}{n}f_n(g_{n-1}(t_n))}{\sqrt{P_n^2(\mathbf{w}; g_{n-1}(t_n))-\frac{4b_na_{nn}}{n\lambda}}}\frac{\ud g_{n-1}(t_n)}{\ud t_n}\nn\\
    &>& 1-\frac{\frac{(n-1)a_{nn-1}}{n}f_n(g_{n-1}(t_n))}{\frac{(n-1)a_{nn-1}}{n}g_{n-1}(t_n)}\frac{g_{n-1}(t_n)}{t_n}\nn\\
    &=&1-\frac{f_n(g_{n-1}(t_n))}{t_n}=\frac{F(t_n)}{t_n}.\label{rr13}
 \eer

Therefore, we  conclude  from \eqref{rr13} that   $F(\cdot)$ is strictly increasing over $(0, \infty)$, which implies the uniqueness
of the zero of $F(\cdot)$.  Then the proof of Proposition \ref{prop1} is complete. $\square$

\subsection{Constrained minimization}
By  Proposition \ref{prop1},   for any  $\mathbf{w}\in\mathcal{A}$,   we see that    the equations  \eqref{s19}--\eqref{s21}  with respect to $ (c_1, \dots,  c_n)$
 admit a solution   $\big(c_1(\mathbf{w}), \dots, c_n(\mathbf{w})\big)$ given by
 \eqref{s30}--\eqref{s32},  such that  $\mathbf{v}$ defined by
 \be v_i=w_i+c_i(\mathbf{w}),  \quad i=1, \dots,  n,\ee
 satisfies  \eqref{s11a}.

Therefore,  to  find  the critical points of  the functional $I$,  we  consider the functional
  \be
   J(\mathbf{w})\equiv I\big(w_1+c_1(\mathbf{w}), \dots, w_n+c_n(\mathbf{w})\big)=I(\mathbf{w}+\mathbf{c}(\mathbf{w})),  \label{c24}
  \ee
where $\mathbf{w}\in \mathcal{A}.$

Noting that \eqref{s11a} is equivalent to \eqref{as1}, and multiplying \eqref{as1} by $\mathbf{1}^\tau$, we see that
  \ber
    \ito   \mathbf{U}^\tau\tilde{S}(\mathbf{U}-\mathbf{1})\ud x+\frac{\mathbf{1}^\tau\mathbf{b}}{\lambda}=0.\label{s38}
  \eer
  Then   in view of $\tilde{K}^{-1}\mathbf{1}=\mathbf{1}$, that is $\tilde{S}^{-1}\tilde{P}^{-1}\mathbf{1}=\mathbf{1}$, we have
    \ber
     -\ito \mathbf{1}^\tau\tilde{P}^{-1}\left(\mathbf{1}-\mathbf{U}\right)\ud x+\ito (\mathbf{U}-\mathbf{1})^\tau\tilde{S}(\mathbf{U}-\mathbf{1})\ud x+\frac{\mathbf{1}^\tau\mathbf{b}}{\lambda}=0.\label{s39}
  \eer

 Hence from \eqref{s39} we  may write the functional $J$ as
 \ber
  J(\mathbf{w})&=&\frac12\sum\limits_{i=1}^2\ito \partial_i \mathbf{w}^\tau M\partial_i\mathbf{w}\ud x+\frac\lambda2 \ito \mathbf{1}^\tau\tilde{P}^{-1}\left(\mathbf{1}-\mathbf{U}\right)\ud x+\mathbf{b}^\tau \mathbf{c}-\frac{\mathbf{1}^\tau\mathbf{b}}{2}, \label{s40}
 \eer
 which  is  \ber
  J(\mathbf{w})&=&\frac12\sum\limits_{i=1}^2\ito \partial_i \mathbf{w}^\tau M\partial_i \mathbf{w}\ud x+\lambda \sum\limits_{i=1}^n\frac{i(n+1-i)}{2n}\ito\left(1-\re^{c_i}\re^{u_i^0+w_i}\right)\ud x\nn\\
  &&+\sum\limits_{i=1}^nb_ic_i- \frac12\sum\limits_{i=1}^nb_i. \label{s40}
 \eer

 It is easy to  see that the functional $J$ is Frech\'{e}t differentiable in the interior of $\mathcal{A}$.  If we find a  minimizer  $\mathbf{w}$ of
$J$, which lies in the interior of $\mathcal{A}$,  then $(\mathbf{w}+\mathbf{c}(\mathbf{w}))$ is  a critical point of $I$.   Hence,  we just
need to  find a minimizer of $J$  in the interior of $\mathcal{A}$, denoted by $\mbox{int}\mathcal{A}$.

 Below we first aim to   find a minimizer of $J$ in $\mathcal{A}$.
We begin by establishing  the following lemma.
  \begin{lemma}\label{lem1}
   For any $\mathbf{w}\in \mathcal{A}$,  there hold
    \ber
     \re^{c_i}\ito\re^{u_i^0+w_i}\ud x\le |\Omega|, \quad i=1, \dots, n,\label{s41}\\
      \re^{c_i}\le 1, \quad i=1, \dots, n. \label{s42}
    \eer
  \end{lemma}

{\bf Proof. }
  Using \eqref{s30}--\eqref{s32} and the constraints \eqref{s28'}, we have
   \ber
    \re^{c_1}&\le& \frac{1}{a_{11}}\left(\frac{a_1}{n}+\frac{(n-1)a_{12}}{n} \re^{c_2}\right),  \label{s43} \\
     \re^{c_i}&\le& \frac{1}{a_{ii}}\left(\frac{a_i}{i(n+1-i)}+ \frac{(i-1)(n+2-i)a_{ii-1}}{2i(n+1-i)}\re^{c_{i-1}}\right.\nn
     \\&&\left.+ \frac{(i+1)(n-i)a_{ii+1}}{2i(n+1-i)}\re^{c_{i+1}}\right), \quad i=2, \dots, n-1,\label{s44}\\
\re^{c_n}&\le&  \frac{1}{a_{nn}}\left( \frac{a_n}{n}+\frac{(n-1)a_{nn-1}}{n}\re^{c_{n-1}}\right). \label{s45}
   \eer

 Let us define an $n\times n$ tridiagonal  matrix $A$  by
 \ber
 A\equiv
\begin{pmatrix}
 na_{11}&-(n-1)a_{12}&0 &\cdots&0\\
 -\frac n2a_{21}&2(n-1)a_{22}&-\frac{3(n-2)}{2}a_{23}&\dots&0\\
  \quad \ddots&\ddots&\ddots&\quad &\vdots\\
   \dots&-A_{ii-1} &A_{ii}&-A_{ii+1}&\cdots\\
\vdots&\ddots&\ddots&\ddots&\vdots\\
0&\dots&-\frac{3(n-2)}{2}a_{n-1n-2}&2(n-1)a_{n-1n-1}&-\frac n2a_{n-1n}\\
 0& \cdots&0& -(n-1)a_{nn-1}&na_{nn}
\end{pmatrix},\label{s45a}
\eer

where
 \ber
&& A_{ii-1}\equiv  \frac{(i-1)(n+2-i)}{2}a_{ii-1}, \quad A_{ii}\equiv i(n+1-i) a_{ii},\\
  &&A_{ii+1}\equiv \frac{(i+1)(n-i)}{2}a_{ii+1}, \quad i=3,\dots, n-2,
 \eer
and we use the notation \eqref{s24b}.

We use the convention that for any  two vectors ${\bf a}=(a_1, \dots, a_n)^\tau$ and ${\bf b}=(b_1, \dots, b_n)^\tau$ we write ${\bf a}\le {\bf b}$ if $a_i\le b_i$ for all $i=1, \dots, n$.

Thus the inequalities \eqref{s43}--\eqref{s45} can be rewritten as
 \ber
  A(\re^{c_1},\dots, \re^{c_n})^\tau\le (a_1, \dots, a_n)^\tau. \label{s45b}
 \eer

Now we denote the adjugate  matrices of $A$ and $\tilde{K}$  by
 \ber
  A^*&=&(A^*_{ij})_{n\times n},\\
  \tilde{K}^*&=&(\tilde{K}^*_{ij})_{n\times n},
 \eer
respectively.

  By using H\"{o}lder's  inequality and an  induction argument we see that all the entries of $A^*$ and the determinant of $A$ are positive.
  Then we may express $A^{-1}$ as
   \ber
    A^{-1}=\frac1{{\det A}}{A^*},
   \eer
   whose entries are all positive.
   Hence  from \eqref{s45b} we obtain
    \berr
 \begin{pmatrix}
  \re^{c_1}\\ \vdots\\ \re^{c_i}\\\vdots\\\re^{c_n}
\end{pmatrix}
 \le A^{-1} \begin{pmatrix}
 a_1\\ \vdots\\ a_i\\\vdots\\a_n
\end{pmatrix}
\eerr
\ber
=\begin{pmatrix}
  \frac{\sum\limits_{j=1}^na_jA_{1j}^*}{na_{11}A_{11}^*-(n-1)a_{12}A_{21}^*}\\
 \vdots\\  \frac{\sum\limits_{j=1}^na_jA_{ij}^*}{-\frac{(i-1)(n+2-i)}{2}a_{ii-1}A_{i-1i}^*+ i(n+1-i)a_{ii}A_{ii}^*-\frac{(i+1)(n-i)}{2}a_{ii+1}A_{i+1i}^*}\\
 \vdots\\
 \frac{\sum\limits_{j=1}^na_jA_{nj}^*}{-(n-1)a_{nn-1}A_{n-1n}^*+na_{nn}A_{nn}^*}
\end{pmatrix}.
    \eer
Therefore,  by repeatedly using H\"{o}lder's inequality,  we have
        \ber
 \begin{pmatrix}
  \re^{c_1}a_1\\ \vdots\\ \re^{c_i}a_i\\\vdots\\\re^{c_n}a_n
\end{pmatrix}
 &\le& \begin{pmatrix}\label{sa1}
  \frac{\sum\limits_{j=1}^na_1a_jA_{1j}^*}{na_{11}A_{11}^*-(n-1)a_{12}A_{21}^*}\\
 \vdots\\  \frac{\sum\limits_{j=1}^na_ia_jA_{ij}^*}{-\frac{(i-1)(n+2-i)}{2}a_{ii-1}A_{i-1i}^*+ i(n+1-i)a_{ii}A_{ii}^*-\frac{(i+1)(n-i)}{2}a_{ii+1}A_{i+1i}^*}\\
 \vdots\\
 \frac{\sum\limits_{j=1}^na_na_jA_{nj}^*}{-(n-1)a_{nn-1}A_{n-1n}^*+na_{nn}A_{nn}^*}
\end{pmatrix}\nn\\
&\le& \begin{pmatrix}
  \frac{|\Omega|\prod\limits_{j=1}^n a_{jj}\sum\limits_{j=1}^n\tilde{K}_{1j}^*}{\prod\limits_{j=1}^n a_{jj}\det \tilde{K}}\\
 \vdots\\  \frac{|\Omega|\prod\limits_{j=1}^n a_{jj}\sum\limits_{j=1}^n\tilde{K}_{ij}^*}{\prod\limits_{j=1}^n a_{jj}\det \tilde{K}}\\
 \vdots\\
  \frac{|\Omega|\prod\limits_{j=1}^n a_{jj}\sum\limits_{j=1}^n\tilde{K}_{nj}^*}{\prod\limits_{j=1}^n a_{jj}\det \tilde{K}}
\end{pmatrix}
= |\Omega|\begin{pmatrix}
 \sum\limits_{j=1}^n(\tilde{K}^{-1})_{1j}\\
 \vdots\\  \sum\limits_{j=1}^n(\tilde{K}^{-1})_{ij}\\
 \vdots\\
  \sum\limits_{j=1}^n(\tilde{K}^{-1})_{nj}
\end{pmatrix}=|\Omega|\tilde{K}^{-1}\mathbf{1}=|\Omega|\mathbf{1}.\label{sa2}
    \eer
   Then the proof of Lemma \ref{lem1} is complete. $\square$

\begin{lemma}\label{lem2}
 For any $\mathbf{w}\in\mathcal{A}$  and  $s\in (0, 1)$, there holds
  \ber
   \ito\re^{u_i^0+w_i}\ud x&\le& \left(\frac{\lambda}{4nb_i}\right)^{\frac{1-s}{s}}\left(\ito\re^{su_i^0+sw_i}\ud x\right)^{\frac1s},\quad i=1, \dots,  n.\label{s46}
  \eer
\end{lemma}

For a proof of this lemma, see \cite{nota, nota2}.

To proceed further,   we  need the well-known Moser--Trudinger inequality
(see\cite{aubi,font})
 \be
 \ito \re^{w}\ud x \le C\exp\left(\frac{1}{16\pi}\|\nabla w\|_2^2\right), \quad \forall\, w\in \dot{W}^{1, 2}(\Omega),\label{c32}
 \ee
 where $C$ is a positive constant depending on $\Omega$ only.

Noting  that the matrices $M$ and $\tilde{S}$,  defined by  \eqref{s5} and  \eqref{s9}, are both positive definite,   we have the following  coercive estimate for $J$.
\begin{lemma}\label{lem3}
  For any $\mathbf{w}\in \mathcal{A}$  there exists a positive constant    $C$ independent of $\lambda$  such that
  \ber
   J(\mathbf{w})\ge  \frac{\alpha_0}{4}\sum\limits_{i=1}^n\|\nabla w_i\|_2^2-C(\ln\lambda+1),  \label{s47}
  \eer
  where $\alpha_0>0$ is the smallest eigenvalue of $M$.
\end{lemma}

{\bf Proof.}   Since the matrices $M$ and $\tilde{S}$,  defined by  \eqref{s5} and  \eqref{s9}, are both positive definite, denoting by $\alpha_0$ the smallest eigenvalue
of $M$, we have
 \ber
  J(\mathbf{w})\ge \frac{\alpha_0}{2} \sum\limits_{i=1}^n\|\nabla w_i\|_2^2+\sum\limits_{i=1}^nb_ic_i. \label{ss1}
 \eer

Using \eqref{s30}--\eqref{s32}, we have
\ber
 \re^{c_i}\ge \frac{a_i}{2i(n+1-i)a_{ii}}=\frac{\ito\re^{u_i^0+w_i}\ud x}{2i(n+1-i)\ito\re^{2u_i^0+2w_i}\ud x}.
\eer

Then by the constraints \eqref{s28'}  we obtain
 \ber
 \re^{c_i}\ge \frac{2nb_i}{i(n+1-i)\lambda\ito\re^{u_i^0+w_i}\ud x},
 \eer
 which implies
 \ber
  c_i\ge \ln\frac{2nb_i}{i(n+1-i)}-\ln\lambda-\ln\ito\re^{u_i^0+w_i}\ud x, \quad i=1, \dots, n. \label{ss2}
 \eer

In view of Lemma \ref{lem2}, and the Moser--Trudinger inequality, we estimate the last term in \eqref{ss2} with
 \ber
  &&\ln\ito\re^{u_i^0+w_i}\ud x\nn\\
  &&\le \frac1s\ln\ito\re^{su_i^0+sw_i}\ud x+\frac{1-s}{s}(\ln\lambda-\ln4nb_i)\nn\\
  &&\le\frac1s\left(\ln C+s\max\limits_{\Omega}u_i^0+\frac{s^2}{16\pi}\|\nabla w_i\|_2^2\right)+\frac{1-s}{s}(\ln\lambda-\ln4nb_i)\nn\\
  &&\le\frac{s}{16\pi}\|\nabla w_i\|_2^2+\frac{1-s}{s}(\ln\lambda-\ln4nb_i)+\frac1s \ln C+\max\limits_{\Omega}u_i^0, \quad i=1, \dots, n.\quad\label{ss3}
 \eer
 Then inserting \eqref{ss3} into \eqref{ss2} gives
  \ber
   c_i&\ge&-\frac{s}{16\pi}\|\nabla w_i\|_2^2-\frac1s(\ln\lambda+\ln C-\ln4nb_i)\nn\\
&&\quad -\ln2i(n+1-i)-\max\limits_{\Omega}u_i^0,\quad  i=1, \dots, n.\quad \label{ss4}
  \eer

 Hence combining \eqref{ss1} and \eqref{ss4}  we get
  \ber
   J(\mathbf{w})&\ge&  \left(\frac{\alpha_0}{2}-\frac{s}{16\pi} \right)\sum\limits_{i=1}^n\|\nabla w_i\|_2^2-\frac1s\sum\limits_{i=1}^nb_i(\ln\lambda+\ln C-\ln4nb_i)\nn\\
   &&-\sum\limits_{i=1}^nb_i\left(\ln2i[n+1-i]+\max\limits_{\Omega}u_i^0\right),
  \eer
 which concludes the lemma with taking $s$ suitably small.  $\square$

  Noting that $J$  is weakly lower semicontinuous in $\mathcal{A}$ and    using Lemma \ref{lem3},  we infer that $J$  has a minimizer in $\mathcal{A}$.

 \subsection{Interior minimizer}
    In the sequel  we show that the minimizer of $J$ obtained above  is an interior point of $\mathcal{A}$ when $\lambda$ is
suitably large.  To this end, we  first estimate the value of the functional $J$ on the boundary of $\mathcal{A}$.
 \begin{lemma}\label{lem4}
  On the boundary of $\mathcal{A}$  there exists a constant $C>0$ independent of $\lambda$ such that
   \ber
    \inf\limits_{\mathbf{w}\in\partial\mathcal{A}} J(\mathbf{w}) &\ge& \frac{|\Omega|\lambda}{2}-C(\ln\lambda+\sqrt{\lambda}+1).\label{s48}
   \eer
 \end{lemma}
{\bf Proof.}
  On the boundary of $\mathcal{A}$,  at least one  of the following  $n$ conditions occurs:
  \ber
   a_i^2&=& \frac{4nb_ia_{ii}}{\lambda}, \quad i=1,  \dots, n.\label{s47a}
 \eer

Without loss of generality, if $i=1$ in (\ref{s47a}), then using \eqref{sa1} and H\"{o}lder's inequality,  we conclude
 \ber
  \re^{c_1}a_1&\le& \frac{\sum\limits_{j=1}^na_1a_jA_{1j}^*}{na_{11}A_{11}^*-(n-1)a_{12}A_{21}^*}\nn\\
  &\le& \frac{a_1^2 \prod\limits_{j=2}^2a_{jj}\tilde{K}_{11}^*}{a_{11}\prod\limits_{j=2}^2a_{jj}\det \tilde{K}}+\frac{a_1\sqrt{|\Omega|}\prod\limits_{j=2}^2a_{jj}\sum\limits_{j=2}^n\tilde{K}_{1j}^*}{\sqrt{a_{11}}\prod\limits_{j=2}^2a_{jj}\det \tilde{K}}\nn\\
  &=& (\tilde{K}^{-1})_{11}\frac{a_1^2}{a_{11}}+\sqrt{|\Omega|}\frac{a_1}{\sqrt{a_{11}}}\sum\limits_{j=2}^n(\tilde{K}^{-1})_{1j}\nn\\
  &=&\frac{8nb_1}{(n+1)\lambda}+\frac{2(n-1)\sqrt{nb_1|\Omega|}}{(n+1)\sqrt\lambda}. \label{ss5}
  \eer

Hence using Lemma \ref{lem1} and \eqref{ss5}, we infer that
 \ber
  &&\lambda \sum\limits_{i=1}^n\frac{i(n+1-i)}{2n}\ito\left(1-\re^{c_i}\re^{u_i^0+w_i}\right)\ud x\nn\\
  &&\ge \frac{|\Omega|\lambda}{2}-\frac{1}{n+1}\left(4nb_1+[n-1]\sqrt{nb_1|\Omega|\lambda}\right). \,\label{ss6}
 \eer

For other cases, we may get similar estimates as \eqref{ss6}.

Then, estimating $c_i, i=1, \dots, n$, as done in Lemma \ref{lem3},  we obtain  desired estimate \eqref{s48}.   $\square$

At this point,  we need to  find some  suitable  test   function, which lies in  the interior of $\mathcal{A}$.  We aim to compare the values of the functional at the test function with that on the boundary of $\mathcal{A}$.

 It was proved in \cite{taran96} that   for $\mu>0$ sufficiently large, the problem
 \ber
  \Delta v=\mu\re^{u_i^0+v}(\re^{u_i^0+v}-1)+\frac{4\pi N_i}{|\Omega|} \quad  \text{in}\quad \Omega, \quad i=1, \dots, n, \label{s49}
 \eer
 admits  a solution $v_i^\mu,$ satisfying
  $u_i^0+v_i^\mu<0$ in $\Omega$, $c_i^\mu=\frac{1}{|\Omega|}\ito v_i^\mu\ud x \to 0$, and   $w_i^\mu=v_i^\mu-c_i^\mu\to -u_i^0$ pointwise as $\mu\to \infty$, $i=1, \dots,  n$.
  In particular, we have the limits
 \ber
   \lim\limits_{\mu\to\infty}\ito \re^{u_i^0+w_i^\mu}\ud x= |\Omega|, \quad \lim\limits_{\mu\to\infty}\ito \re^{2u_i^0+2w_i^\mu}\ud x= |\Omega|, \quad i=1, \dots,  n.\label{s51}
 \eer

Let us introduce  an $n\times n$ tridiagonal  matrix $\tilde{A}(\mathbf{w}^\mu)$  defined as
\ber
\tilde{A}(\mathbf{w}^\mu)\equiv
\begin{pmatrix}
 na_{11}&-(n-1)|\Omega|&0 &\cdots&0\\
  -\frac n2|\Omega|&2(n-1)a_{22}&-\frac{3(n-2)}{2}|\Omega|&\dots&0\\
 \ddots&\ddots&\ddots&\vdots&\vdots\\
   \dots&-\frac{(i-1)(n+2-i)}{2}|\Omega| &i(n+1-i)a_{ii}&-\frac{(i+1)(n-i)}{2}|\Omega|&\cdots\\
\vdots&\ddots&\quad\ddots&\ddots&\vdots\\
0&\dots&-\frac{3(n-2)}{2}|\Omega|&2(n-1)a_{n-1n-1}&-\frac n2|\Omega|\\
 0& \cdots&0& -(n-1)|\Omega|&na_{nn}
\end{pmatrix},\label{ss7}
\eer
where we use the notation \eqref{s24b} with $a_{ii}=a_{ii}(w_i^\mu), i=1, \dots, n.$
Then  we see from \eqref{s51} that
 \ber
  \lim\limits_{\mu\to\infty} \tilde{A}(\mathbf{w}^\mu)=|\Omega|\tilde{K}. \label{ss8}
 \eer

 It is easy to see from Jensen's inequality that all the cofactors and the determinant of  $\tilde{A}(\mathbf{w}^\mu)$ are positive, which implies in particular $\tilde{A}(\mathbf{w}^\mu)$ is invertible and all the
 entries of $\tilde{A}^{-1}(\mathbf{w}^\mu)$ are positive.
 Then it follows from \eqref{ss8} that
  \ber
  \lim\limits_{\mu\to\infty} \tilde{A}^{-1}(\mathbf{w}^\mu)=\frac{1}{|\Omega|}\tilde{K}^{-1}. \label{ss9}
 \eer

Hence we conclude from  the   limits \eqref{s51}, \eqref{ss9}, and  the   definition of $\mathcal{A}$ that,     for a fixed   $\tilde{\lambda}_0>0$ large and for  any $\vep\in(0, 1)$,   there exists a  $\mu_\vep\gg1$, such that
\be
\mathbf{w}^{\mu_\vep}=(w_1^{\mu_\vep},  \dots, w_n^{\mu_\vep})\in \mbox{int}\mathcal{A},\label{s52}
\ee
 for every $\lambda>\tilde{\lambda}_0$,  and there hold
  \ber
  && {\rm diag}\Big\{a_{11}(w_1^{\mu_\vep}), \dots, a_{nn}(w_n^{\mu_\vep})\Big\}<2|\Omega|{\rm diag}\Big\{1, \dots, 1\Big\},\label{ss10'}\\
   &&\frac{ (1-\vep)}{|\Omega|}\tilde{K}^{-1}\le \tilde{A}^{-1}(\mathbf{w}^{\mu_\vep})\le \frac{ (1+\vep)}{|\Omega|}\tilde{K}^{-1}<\frac{ 2}{|\Omega|}\tilde{K}^{-1}.\label{ss10}
  \eer
  Here we use the convention that for $n\times n$ matrices $G=(g_{ij})_{n\times n}$ and $H=(h_{ij})_{n\times n}$, we write   $G\le H$ if $g_{ij}\le h_{ij}$ for all $i, j=1, \dots, n$.

Now we prove the following comparison result.

\begin{lemma}\label{lem5}
 When $\lambda>0$ is  suitably large,  for the test vector  $\mathbf{w}^{\mu_\vep}$ given by \eqref{s52},     there holds
   \ber
    J(\mathbf{w}^{\mu_\vep})-\inf\limits_{\mathbf{w}\in\partial\mathcal{A}}J(\mathbf{w})<-1.\label{s56}
   \eer
\end{lemma}

{\bf  Proof.}\quad For $\mathbf{w}^{\mu_\vep}$ given by \eqref{s52},  using Proposition \ref{prop1},  we may get a corresponding vector  $(c_1(\mathbf{w}^{\mu_\vep}), \dots,  c_n(\mathbf{w}^{\mu_\vep}))$ defined by \eqref{s30}--\eqref{s32}.
 Then applying  Jensen's inequality      and \eqref{s30}--\eqref{s32},  we obtain
  \ber
  \re^{c_1(\mathbf{w}^{\mu_\vep})}&=&\frac{P_1\big(\mathbf{w}^{\mu_\vep};  \re^{c_2(\mathbf{w}^{\mu_\vep})}\big)}{2a_{11}(w_1^{\mu_\vep})}\left(1+\sqrt{1-\frac{4b_1a_{11}(w_1^{\mu_\vep})}{n\lambda P_1^2\big(\mathbf{w}^{\mu_\vep};  \re^{c_2(\mathbf{w}^{\mu_\vep})}\big)}}\right)\nn\\
  &\ge&\frac{P_1\big(\mathbf{w}^{\mu_\vep};  \re^{c_2(\mathbf{w}^{\mu_\vep})}\big)}{a_{11}(w_1^{\mu_\vep})}-\frac{2b_1}{n\lambda P_1\big(\mathbf{w}^{\mu_\vep};  \re^{c_2(\mathbf{w}^{\mu_\vep})}\big)}\nn\\
   &\ge& \frac{|\Omega|}{a_{11}(w_1^{\mu_\vep})}{\left(\frac1n+\frac{n-1}{n}\re^{c_2(\mathbf{w}^{\mu_\vep})}\right)}-\frac{2b_1}{\lambda|\Omega|}.\label{s53}
  \eer
  Analogously,
  \ber
   \re^{c_i(\mathbf{w}^{\mu_\vep})}&\ge& \frac{|\Omega|}{a_{ii}(w_i^{\mu_\vep})}\left\{\frac{1}{i(n+1-i)}+ \frac{(i-1)(n+2-i)}{2i(n+1-i)}\re^{c_{i-1}(\mathbf{w}^{\mu_\vep})} \right.\nn\\
   &&\left.+ \frac{(i+1)(n-i)}{2i(n+1-i)}\re^{c_{i+1}(\mathbf{w}^{\mu_\vep})}\right\}-\frac{2nb_i}{\lambda i(n+1-i)|\Omega|},\quad i=2, \dots, n-1, \label{s54}\\
   \re^{c_n(\mathbf{w}^{\mu_\vep})}
  &\ge& \frac{|\Omega|}{a_{nn}(w_n^{\mu_\vep})}{\left(\frac1n+\frac{n-1}{n}\re^{c_{n-1}(\mathbf{w}^{\mu_\vep})}\right)}-\frac{2b_n}{\lambda|\Omega|},\label{s55}
 \eer
where  we use the notation \eqref{s24b} with the understanding that $a_{ii}=a_{ii}(w_i^{\mu_\vep}), \,i=1, \dots,  n$.

Noting the definition  \eqref{ss7},  we   may express the inequalities \eqref{s53}--\eqref{s55}  equivalently as
 \ber
 \tilde{A}(\mathbf{w}^{\mu_\vep})\left(\re^{c_1(\mathbf{w}^{\mu_\vep})},   \dots, \re^{c_n(\mathbf{w}^{\mu_\vep})}\right)^\tau
  \ge |\Omega|\mathbf{1}-\frac{2n}{\lambda|\Omega|}{\rm diag}\Big\{a_{11}(w_1^{\mu_\vep}), \dots, a_{nn}(w_n^{\mu_\vep})\Big\}\mathbf{b}.\label{ss12}
 \eer
Since all the  entries of $\tilde{A}^{-1}(\mathbf{w}^{\mu_\vep})$ are positive,   we infer from \eqref{ss12}, \eqref{ss10'}, and \eqref{ss10} that
 \ber
   && \left(\re^{c_1(\mathbf{w}^{\mu_\vep})},   \dots, \re^{c_n(\mathbf{w}^{\mu_\vep})}\right)^\tau \nn\\
    &&\ge \tilde{A}^{-1}(\mathbf{w}^{\mu_\vep})\left( |\Omega|\mathbf{1}-\frac{2n}{\lambda|\Omega|}{\rm diag}\Big\{a_{11}(w_1^{\mu_\vep}), \dots, a_{nn}(w_n^{\mu_\vep})\Big\}\mathbf{b}\right)\nn\\
    &&\ge \tilde{A}^{-1}(\mathbf{w}^{\mu_\vep})\left( |\Omega|\mathbf{1}-\frac{4n}{\lambda}\mathbf{b}\right)\nn\\
 &&\ge  (1-\vep) \tilde{K}^{-1} \mathbf{1}- \frac{ 8n}{\lambda|\Omega|}\tilde{K}^{-1}\mathbf{b}\nn\\
  & &= (1-\vep)\mathbf{1}-\frac{ 8n}{\lambda|\Omega|}R^{-1}K^{-1}\mathbf{b}.\label{ss13}
 \eer
Then, it follows from \eqref{ss13} that
 \ber
 \ito\left(1-\re^{c_i(\mathbf{w}^{\mu_\vep})}\re^{u_i^0+w_i^{\mu_\vep}}\right)\,\ud x&\le&|\Omega|\vep+\frac{16n}{\lambda i(n+1-i)}\sum\limits_{j=1}^n(K^{-1})_{ij}b_j, \quad i=1, \dots, n.\qquad \label{ss14}
 \eer

At this point using \eqref{s41} and \eqref{ss14} we conclude that,    for any small $\vep>0$, there exists a positive constant $C_\vep$ independent of $\lambda$ such
that
\ber
 J(\mathbf{w}^{\mu_\vep})\le \frac{ (n+1)(n+2)|\Omega|\lambda}{12 }\vep+C_\vep. \label{ss15}
\eer
 Consequently, from \eqref{ss15} and  Lemma \ref{lem5}, we  infer   that
  \ber
    J(\mathbf{w}^{\mu_\vep})-\inf\limits_{\mathbf{w}\in\partial\mathcal{A}} J(\mathbf{w})\le\frac{|\Omega|\lambda}{2}\left( \frac{[n+1][n+2]}{6}\vep-1\right)+C(\ln\lambda+\sqrt{\lambda}+1),\label{ss16}
  \eer
where $C>0$ is a  constant independent of $\lambda$.
Hence,  by taking   $\vep$ suitably small and $\lambda$ sufficiently large in \eqref{ss16},  we obtain the desired estimate    \eqref{s56}.   $\square$

 Now we  may infer from Lemma \ref{lem3} and \ref{lem5}   that there exists a  $ \lambda_1>0$ such that, for every   $\lambda>\lambda_1$, the functional $J$  admits a minimizer
  \be
  \mathbf{w}^\lambda\equiv(w_1^\lambda,  \dots, w_n^\lambda)^\tau\in  \mbox{int} \mathcal{A}. \label{ss17}
  \ee
\subsection{Solution to the original system}
Since we use a constrained minimization to get a  minimizer $\mathbf{w}^\lambda$ of $J(\mathbf{w})=I(\mathbf{w}+\mathbf{c}(\mathbf{w}))$ in the subspace of $W^{1, 2}(\Omega)$,  it is  not obvious that whether  $\mathbf{w}^\lambda$ gives rise to a  solution of the system \eqref{s8}.
Here we show that   $\mathbf{v}^\lambda\equiv(v_1^\lambda, \dots, v_n^\lambda)$ defined by
   \ber
    v_i^\lambda=w_i^\lambda+c_i(\mathbf{w}^\lambda), \quad i=1, \dots,  n,\label{s57}
   \eer
is  actually  a  solution of the system \eqref{s8}.
\begin{lemma}\label{lem6}
 Let $\mathbf{w}$ be a  minimizer of $J$ in $\mbox{\rm int}\mathcal{A}$ and the corresponding vector $\mathbf{c}(\mathbf{w})$  be determined by \eqref{s30}--\eqref{s32}. Then
\ber
\mathbf{v}=\mathbf{c}(\mathbf{w})+\mathbf{w}
\eer
must be  a solution of the system  \eqref{s8}.
\end{lemma}

{\bf Proof.} \,   Since   $\mathbf{w}$  is an interior minimizer of $J$ in $\mathcal{A}$,  the Fr\'{e}chet derivative of $J(\mathbf{w})=I(\mathbf{w}+\mathbf{c}(\mathbf{w}))$ at $\mathbf{w}$
 should  be zero,
 \ber
  [{\rm d }I(\mathbf{w}+\mathbf{c}(\mathbf{w}))]\mathbf{f}=0 \quad  \text{for any }\,   \mathbf{f}\in \dot{W}^{1,2}(\Omega).\label{ss21}
 \eer
 By the expression of $I$ \eqref{s11}, we rewrite \eqref{ss21} in an explicit form
 \ber
  &&\ito\left(\sum\limits_{i=1}^2\partial_i\mathbf{f}^\tau M\partial_i\mathbf{w}+\lambda\mathbf{f}^\tau\mathrm{U}\tilde{S}[\mathbf{U}-\mathbf{1}]\right)\ud x\nn\\
  && +[D_{\mathbf{f}}\mathbf{c}(\mathbf{w})]^\tau\ito\left( \lambda\mathrm{U}\tilde{S}[\mathbf{U}-\mathbf{1}]+\frac{\mathbf{b}}{|\Omega|}\right)\ud x=0   \quad \text{for any }\, \mathbf{f}\in \dot{W}^{1,2}(\Omega), \label{ss22}
 \eer
where
 \ber
 D_{\mathbf{f}}\mathbf{c}(\mathbf{w})=\frac{\ud}{\ud t}\mathbf{c}(\mathbf{w}+t\mathbf{f})|_{t=0},
 \eer
is the directional derivative of $\mathbf{c}$ at $\mathbf{w}$ along the direction $\mathbf{f}$, and the notation \eqref{s1'} is used.

Then  we  use \eqref{as1} to  reduce   \eqref{ss22} into
 \be
  \ito\left(\sum\limits_{i=1}^2\partial_i\mathbf{f}^\tau M\partial_i\mathbf{w}+\lambda\mathbf{f}^\tau\mathrm{U}\tilde{S}[\mathbf{U}-\mathbf{1}]\right)\ud x=0.\label{ss23}
 \ee

 Denote by  $L^2(\Omega)$  the scalar-valued or $n$-vector-valued function space of $\Omega$-period   $L^2$-functions and decompose  $L^2(\Omega)$ as
\be L^2(\Omega)=\mathbb{R}^n\oplus Y,\ee
where
\be
 Y=\left\{\mathbf{f}\Big|\, \mathbf{f}\in L^2(\Omega), \quad \ito \mathbf{f}\ud x=\mathbf{0}\right\}.
\ee

We  select   a vector $\mathbf{d}\in \mathbb{R}^n$ such that
\be
\lambda\mathrm{U}\tilde{S}(\mathbf{U}-\mathbf{1})+\mathbf{d}\in Y.
\ee
 Hence the relation  $\dot{W}^{1,2}(\Omega)\subset Y$ and \eqref{ss23} lead to
  \ber
  0&=&\ito\left\{\sum\limits_{i=1}^2\partial_i\mathbf{f}^\tau M\partial_i\mathbf{w}+\mathbf{f}^\tau\left(\lambda\mathrm{U}\tilde{S}[\mathbf{U}-\mathbf{1}]+\mathbf{d}\right)\right\}\ud x\nn\\
   &=& \ito\left\{\sum\limits_{i=1}^2\partial_i(\mathbf{f}+\mathbf{g})^\tau M\partial_i\mathbf{w}+(\mathbf{f}+\mathbf{g})^\tau\left(\lambda\mathrm{U}\tilde{S}[\mathbf{U}-\mathbf{1}]+\mathbf{d}\right)\right\}\ud x,\label{ss24}
 \eer
for any $\mathbf{g}\in \mathbb{R}^n$.  Consequently,  we have
 \be
  \ito\left\{\sum\limits_{i=1}^2\partial_i\mathbf{h}^\tau M\partial_i\mathbf{w}+\mathbf{h}^\tau\left(\lambda\mathrm{U}\tilde{S}[\mathbf{U}-\mathbf{1}]+\mathbf{d}\right)\right\}\ud x=0  \quad \text{for any }\, \mathbf{h}\in W^{1,2}(\Omega). \label{ss25}
 \ee
  Then we conclude from \eqref{ss25} that  $\mathbf{w}$ is a smooth solution of the system
 \be
 \Delta M\mathbf{w}=\lambda\mathrm{U}\tilde{S}(\mathbf{U}-\mathbf{1})+\mathbf{d},\label{ss26}
 \ee
 which, after being integrated over $\Omega$,   gives us
  \be
 \lambda \ito\mathrm{U}\tilde{S}(\mathbf{U}-\mathbf{1})\ud x+\mathbf{d}|\Omega|=\mathbf{0}.\label{ss27}
 \ee
   Hence we infer from \eqref{ss27} and  \eqref{as1}  that
   \be
    \mathbf{d}=\frac{\mathbf{b}}{|\Omega|}. \label{ss28}
   \ee
 Combining  \eqref{ss26} and \eqref{ss28}  we see   that  $\mathbf{v}=\mathbf{c}(\mathbf{w})+\mathbf{w}$ is a solution of \eqref{s8}.  Thus the  proof of Lemma \ref{lem6} is complete. $\square$

 At this stage we infer    from Lemma \ref{lem6} that  when $\lambda>\lambda_1$, $\mathbf{v}^\lambda$ defined by \eqref{s57} is a solution of \eqref{s8}.  Therefore   part (ii) of Theorem \ref{th2}
follows.

\subsection{Asymptotic behavior and quantized integrals}

In this subsection we study the asymptotic behavior of the solution $\mathbf{v}^\lambda$ of \eqref{s8} defined by \eqref{s57}  when $\lambda\to \infty$  and
establish the quantized integrals as stated in Theorem \ref{th2}.

 \begin{lemma}\label{lem7}
  Let $\mathbf{v}^\lambda$  be the solution of  \eqref{s8}  given by  \eqref{s57}.    Then there holds
   \be
    \lim\limits_{\lambda\to\infty}\ito\left(\re^{u_i^0+v_i^\lambda}-1\right)^2\ud x=0, \quad i=1, \dots, n. \label{ss17a}
   \ee
 \end{lemma}
{\bf Proof. }\quad
   Since $J$ achieves its minimum at $\mathbf{w}^\lambda\in  \mbox{int}\mathcal{A}$,   we  see from \eqref{ss15} that, for any $\vep\in (0, 1)$, there exist   constants   $\lambda_\vep>0$ and $C_\vep>0$  such that
   \ber
    J(\mathbf{w}^\lambda)=\inf\limits_{\mathbf{w}\in\mathcal{A}}J(\mathbf{w})\le \frac{ (n+1)(n+2)\lambda|\Omega|}{12 }\vep+C_\vep \quad  \text{for all}\quad  \lambda>\lambda_\vep. \label{ss18}
   \eer
 Noting the matrix $\tilde{S}$ (defined by \eqref{s5}) is positive definite  and denoting   the smallest eigenvalue of  $\tilde{S}$ by    $\beta_0>0$, we have
\ber
 \ito(\mathbf{U}-\mathbf{1})^\tau\tilde{S}(\mathbf{U}-\mathbf{1})\ud x\ge \beta_0\sum\limits_{i=1}^n\ito\left(\re^{u_i^0+v_i}-1\right)^2\ud x.\label{ss19}
\eer
Therefore, in view of \eqref{s39}, \eqref{s40}, and  \eqref{ss19},  and estimating  $c_i^\lambda$ as  that  in  Lemma \ref{lem4},  we conclude that
 \ber
 J(\mathbf{w}^\lambda)\ge \frac{\beta_0\lambda}{2}\sum\limits_{i=1}^n\ito\left(\re^{u_i^0+v_i^\lambda}-1\right)^2\ud x-C(\ln\lambda+1),\label{ss20}
\eer
where  $C>0$ is a constant independent of $\lambda$.

Then combining  \eqref{ss18} and \eqref{ss20} leads to
\berr
 \limsup\limits_{\lambda\to \infty}\sum\limits_{i=1}^n\ito\left(\re^{u_i^0+v_i^\lambda}-1\right)^2\ud x\le \frac{(n+1)(n+2)|\Omega|}{6\beta_0}\vep\quad \text{for any}\quad \vep>0,
\eerr
which  implies the desired conclusion \eqref{ss17a}. The proof of Lemma \ref{lem7} is complete. $\square$

Hence part (iii) of Theorem \ref{th2} follows from  Lemma \ref{lem7}.

To establish  the quantized integrals \eqref{t11}, we  just need to    integrate the   equations    \eqref{s2a} over $\Omega$.

  The proof of Theorem \ref{th2} is complete.

\section{Concluding remarks}
\setcounter{equation}{0}

We note that the method to establish Theorem \ref{th2}  can be applied to prove an existence theorem for the problem   \eqref{t4a} or  \eqref{t5}  when the matrix $\tilde{K}$  assumes  a more general tridiagonal matrix $\hat{K}$
 form,
  \ber
\hat{K}\equiv\begin{pmatrix}
 1+\alpha_{12}&-\alpha_{12}&\cdots&\cdots&0\\
 \ddots&\ddots&\ddots&\vdots&\vdots\\
   \cdots&-\alpha_{ii-1}& 1+\alpha_{ii-1}+\alpha_{ii+1}& -\alpha_{ii+1} &\cdots\\
 \vdots&\ddots&\ddots&\ddots&\vdots\\
 0& \cdots&\cdots&-\alpha_{nn-1}&1+\alpha_{nn-1}
\end{pmatrix}, \label{tn1}\eer
where $\alpha_{12}, \alpha_{ii-1}, \alpha_{ii+1},  \alpha_{nn-1}>0, i=2, \dots, n-1$.

 In fact, it is ready to check that all entries of $\hat{K}^{-1}$ are positive, $\hat{K}^{-1}$ satisfies $\hat{K}^{-1}\mathbf{1}=\mathbf{1}$, and
$\hat{K}$ can be decomposed as
 \be
   \hat{K}=\hat{P}\hat{S},
 \ee
where $\hat{P}$ is a diagonal matrix with positive diagonal entries  and $\hat{S}$ is a positive definite matrix.

In particular, when
\ber
 \alpha_{12}&=&n-1=\alpha_{nn-1}, \nn\\
 \alpha_{ii-1}&=&\frac{(i-1)(n+2-i)}{2}, \\
\alpha_{ii+1}&=&\frac{(i+1)(n-i)}{2},\quad   i=2, \dots, n-1,\nn
 \eer
the matrix $\hat{K}$ reduces to $\tilde{K}$ given by \eqref{t7}.

Since the corresponding existence  result  can be stated in a similar formulation as that of Theorem \ref{th2}, the details are omitted here.

It will be of future interest to develop an existence theory when the Cartan matrix is not tridiagonal.
\medskip

\small{
Han was  supported in part by the National Natural Science Foundation of China under grant 11201118  and by the Key Foundation for Henan Colleges under grant 15A110013. Both  authors were supported in part by the National Natural Science Foundation of China under grants 11471100 and 11471099.  
}

\small{

}
\end{document}